\input amstex 
\documentstyle{amsppt}
\nologo
\magnification=1200

\input xy
\xyoption{all}
\CompileMatrices
\loadbold

\define\fib{\twoheadrightarrow}
\define\iso{\tilde\rightarrow}

\define\g{\frak{g}}
\define\zahl{\Bbb Z}

\define\comm{\operatorname{((Comm))}}

\define\poiss{\operatorname{((Poiss))}}

\define\fpoiss{\operatorname{Poiss}}
\topmatter  

\title The structure of smooth algebras in Kapranov's framework for
noncommutative geometry.
\endtitle
\author by Guillermo Corti\~ nas*\endauthor
\affil Departamento de Matem\'atica\\
Facultad de Cs. Exactas y Naturales\\
    Universidad de Buenos Aires\\
       Argentina\endaffil
\thanks (*) Partially supported by grants BID802/OC-AR-PICT 2260 
and UBACyT TW79. This research was carried out while visiting the 
universities of Bielefeld and M\"unster with a Humboldt fellowship.
\endthanks
\address Departamento de
Matem\' atica, 
Ciudad Universitaria Pabell\'on 1,
(1428) Buenos Aires, Argentina.
\endaddress
\email gcorti\@dm.uba.ar\endemail
\leftheadtext{Guillermo Corti\~ nas}
\rightheadtext{Smooth algebras in noncommutative geometry}
\abstract In Kapranov, M. {\it Noncommutative geometry based on commutator
expansions,} J. reine angew. Math {\bf 505} (1998) 73-118, 
 a theory of noncommutative algebraic varieties
was proposed. Here we prove a structure theorem for the noncommutative 
coordinate rings of affine open subsets of such of those varieties which 
are smooth (Theorem 3.4). The theorem describes the local ring of a point 
as a truncation of a quantization of the enveloping Poisson algebra of a 
smooth commutative local algebra. An 
explicit descripition of this quantization is given in Theorem 2.5. 
A description of the $A$- module structure of the Poisson envelope of 
a smooth commutative algebra $A$ was given 
in {\it loc. cit.}, Theorem 4.1.3. However the proof given in {\it loc. cit.}
has a gap. We fix this gap for $A$ local (Theorem 1.4) and prove a weaker 
global result (Theorem 1.6).
\endabstract

\endtopmatter  

\NoBlackBoxes
\document
\head{0. Introduction}\endhead 
We consider associative, unital algebras over a fixed algebraically closed 
field $k$ of characteristic zero. If $R$ is an algebra, then the {\it 
commutator filtration} of $R$ is defined as 
$$
F_0R=R,\qquad F_{n+1}R:=\sum_{p=1}^nF_pRF_{n+1-p}R+\sum_{p=0}^n
<[F_pR,F_{n-p}R]>\tag{1}
$$
(One checks that this filtration is the same as that called $NC$-filtration in [2].)
By definition, $FR$ is the smallest of all descending filtrations $\Cal G$
with $\Cal G_p\Cal G_q\subset\Cal G_{p+q}$, 
$[\Cal G_p,\Cal G_q]\subset \Cal G_{p+q+1}$. An algebra $R$ is 
{\it nilcommutative of order} $d$ if $F_{d+1}R=0$. Thus the nilcommutative
algebras of order $0$ are just the commutative algebras.
We write $NC_d$ for the
category of nilcommutative algebras of order $d$ and algebra homomorphisms
and set $NC=\cup_{d=0}^\infty NC_d$. An algebra $R$ is called {\it $d$-formally smooth}
if $R\in NC_d$ and if $\hom_{NC_d}(R,)$ maps surjections with nilpotent kernel
to surjections, and is {\it $d$-smooth} if it is $d$-formally smooth and
if the commutative algebra $A=R/F_1R$ is essentially of finite type. For
example a $0$-smooth algebra is the same thing as a smooth commutative algebra. For the
remainder of this section $A$ will be a fixed smooth commutative algebra.
It was shown in [2],1.6.1 that there exists a tower of surjective 
homomorphisms:
$$
\dots \fib R_{d+1}\fib R_d\fib \dots\fib R_2\fib R_1\fib R_0=A
$$
such that $R_d$ is $d$-smooth and $R_d/F_dR=R_{d-1}$ ($d\ge 1$). Moreover
it is shown in {\it loc. cit.} that such a tower is unique up to noncanonical
isomorphism. Kapranov further develops a theory of nilcommutative $d$-smooth 
algebraic  varieties based on this `affine' construction.
In this paper we focus on the affine part of Kapranov's work. We
study the structure of the algebras $R_d$ and of their associated graded
Poisson algebras $GR_d:=\oplus_{n=0}^dF_nR_d/F_{n+1}R_d$. A characterization
of $GR_d$ was given in [2] 4.2.1. It was shown that there is an isomorphsim
$$
PA/P_{>d}A\iso GR_d\tag{2}
$$
Here $P:\comm@>>>\poiss$ is left adjoint to the forgetful functor which 
associates to a Poisson algebra its underlying commutative algebra. The algebra
$PA$ turns out to be graded, and $P_{>d}A=\oplus_{n>d}P_nA$. The map 
\thetag{2}
is canonical, and comes from the adjointness property of $P$; if $R\in NC_d$ 
is any algebra with $G_0R=R/F_1R=A$, then the adjunction map $PA@>>>GR$
induced by $A\cong G_0R$ factors through $PA/P_{>d}A$ obtaining \thetag{2}.
Here we prove a converse of Kapranov's result.
We show that if $R_d\in NC_d$ is any algebra with $R_d/F_1R_d=A$ and
such that the canonical map \thetag{2} is an isomorphism, then $R_d$ is 
$d$-smooth (Theorem 3.4).  This means that to give a $d$-smooth algebra
$R$ with $R/F_1R=A$ is the same thing as to give an associative multiplication
$$
\phi=\sum_{r=0}^d\phi_r:PA/P_{>d}A\otimes P/P_{>d}A@>>>PA/P_{>d}A\tag{3}
$$
with $\phi_r$ homogeneous of degree $r$. For $(A,\Cal M)$ local, we give
(Theorem 2.5) a canonical construction which produces an associative product 
$$
B^X(\hbar)=\sum_{p=0}^\infty B^X_p\hbar^p:PA\otimes PA[[\hbar]]@>>>PA[[\hbar]]
\tag{4}
$$
for each regular system of parameters $X=\{x_1,\dots, x_n\}\subset \Cal M$, 
with $B^X_p$ homogeneous of degree $p$ and a bidifferential operator of
order $\le p$. It turns out that, modulo $P_{>d}A$, the evaluated series 
$B^X(1)$ is a finite sum, and gives a product $\phi$ satisfying the 
requirements of \thetag{3}. We
use this product to give a local characterization of $R_d$ (Theorem 3.4). The 
construction of the product \thetag{4} uses a local isomorphism of 
$A$-modules
$$
PA\cong S^AL_+^A\Omega^1_A\qquad(n\ge 0)\tag{5}
$$ 
Between $PA$ and the $A$-symmetric algebra of the Lie subalgebra 
$$
L_+^A\Omega^1_A=[L^A\Omega^1_A,L^A\Omega^1_A]
$$
of the free $A$-Lie algebra $L^A\Omega^1_A$ generated by the module of 
K\"ahler differentials. Theorem 4.1.3 of [2] states that
there is a global isomorphism as that of \thetag{5}; there is however
a gap in the proof. The gap is explained in section 1 below, where it
is also shown how it is fixed for $A$ local (Theorem 1.4). I do not know 
whether
\thetag{5} still holds globally. A weaker version of \thetag{5}
which holds globally is proved in Theorem 1.6; it establishes that 
$PA$ carries a filtration such that the associated graded module is
(globally) isomorphic to the right hand side of \thetag{5}.

When I explained to Kapranov the gap in his proof of \thetag{5}, and
told him the gap could be fixed locally, he suggested that a weaker version
along the lines of that presented here (Theorem 1.6) should hold globally. 
I am thankful to him for this suggestion.

The remainder of this paper is organized as follows. In section 1 we recall
in some detail the construction of the Poisson algebra $PA$, --which we call
the Poisson envelope of $A$ -- explain the
gap in Kapranov's proof of \thetag{5}, and prove it in the local 
case (Theorem 1.4). The section ends with the  
weaker version of \thetag{5} which holds globally (Theorem 1.6). 
Section 2 is devoted to the construction
of the product \thetag{4} (Theorem 2.5). The results of this section can
be seen as the generalization to general local smooth algebras of those 
obtained by Kapranov for localizations of polynomial rings ([2]\S3). Our
approach is however different from that of [2]. In {\it loc. cit.}
the Feynmann-Maslov
calculus was used to describe the product of elements in the tensor algebra
on a finite dimensional vectorspace $V$
in terms of a specific ordered basis. Instead we use the coordinate free
approach of [1], where explicit formulas for this product were obtained for
not necessarily finite dimensional vectorspaces $V$. The same formulas 
apply to the quantized product \thetag{4}. 
In section 3 we prove (Theorem 3.4) that, for $B^X$ as in \thetag{4},
an algebra $R\in NC_d$ is (i)
$d$-smooth $\Leftrightarrow$ (ii) $A:=R/F_1R$ is smooth commutative and 
\thetag{2} is an isomorphism $\Leftrightarrow$(iii) $R$ is locally isomorphic 
to $(PA/P_{>d}A,B^X(1))$ for some regular system of parameters $X$. Part 
(1)$\Rightarrow$(2) of this was proven by Kapranov in [2] 4.2.1; we give a 
new proof.
\bigskip
\head{1. The Poisson envelope of a commutative algebra}\endhead
\bigskip
\subhead{1.0. Two gradings in the symmetric algebra of a free Lie 
algebra}\endsubhead 
If $V$ is a vectorspace, we write $TV$ for the tensor algebra and 
$LV\subset TV$ for the Lie subalgebra it generates. For $V=\oplus_{x\in X}kx$
--the free vectorspace on a set $X$--
$LV$ is the free Lie algebra on $X$. The symmetric algebra $S\g$ of
any Lie algebra $\g$ is viewed as a Poisson algebra via the Poisson bracket
$\{,\}$ induced by the Lie bracket $[,]$ of $\g$. For example
$$
\fpoiss V:=SLV
$$ 
is a free Poisson algebra. 
Fix a vectorspace $V$ and set $L=LV$. We have $L=\oplus_{n\ge 0}L_n$, where
$$
L_0=V,\quad L_{n+1}=[L_0,L_n]\quad (n\ge 0)\tag{6}
$$
Note our grading is the usual one --as defined for example in [3] LA, 
Ch. IV -- shifted down one degree.
Put
$$
|l|_*=n\ \ \text{if}\ \ l\in L_n\quad (n\ge 0)\tag{7}
$$
This grading induces one in the symmetric algebra $S=SL$; we write $S_n$ 
for its homogeneous part of degree $n$. Note that
$S_n$ is not the same thing as the $n$-th symmetric power $S^n=S^nL$. The
latter is the homogenous part of degree $n$ with respect to a different
grading, namely that given by
$$
|l|^*=1\quad\text{if\ \ }l\in L\tag{8}
$$
Put 
$$
L_+=\oplus_{n\ge 1}L_n
$$
We have $S_0L=SV$ and for $n\ge 1$
$$
\gather
S_nL=SV\otimes S_nL_+\\
\\
S_nL_+=\bigoplus_{r\ge 1}\bigoplus_{\gathered
0<i_1<\dots <i_r\\
p_1i_1+\dots +p_ri_r=n\\
p_1,\dots,p_r>0\endgathered} S^{{p}_1}L_{{i}_1}\otimes\dots\otimes
S^{{p}_r}L_{{i}_r}\tag{9}
\endgather
$$
\bigskip
\subhead{1.1. Poisson ideals}\endsubhead A {\it Poisson ideal} in a Poisson
algebra $P$ is a subspace $I\subset P$ which is an ideal for both
the associative and the Lie algebra structures. If $Y\subset P$ is a 
subset, then we put $<Y>$ and $<<Y>>$ for the smallest ideal and
the smallest Poisson ideal containing $Y$. By definition $<Y>\subset <<Y>>$.
In fact $<<Y>>$ is generated as an ideal by the elements of $Y$ and by those
of the form
$$
\{a_1,\{a_2,\dots,\{a_n,y\}\dots\}\}\quad n\ge 1,\ \ a_i\in P,\ \  y\in Y,
\ \ (1\le i\le n)
$$
Furthermore, one checks that if $X\subset P$ generates $P$ as a 
Poisson algebra then for 
$$
g_i(x_1,\dots,x_n; y):=\{x_1,\{x_2,\dots,\{x_i,\{y,\{x_{i+1},\dots,\{x_{n-1},
x_n\}\dots\}\}\}\dots\}\}\tag{10}
$$
we have
$$
<<Y>>=<Y\cup\bigcup_{n=1}^\infty\{g_i(x_1,\dots,x_n;y)\ \ :\ \ 0\le i\le n,
\ \ x_i\in X,\ \ y\in Y\}>\tag{11}
$$
\bigskip
\subhead{1.2. Poisson envelope}\endsubhead Let $A$ be a commutative algebra,
$SA$ the symmetric algebra on its underlying vectorspace, $SA\fib A$ the 
canonical projection, $IA$ its kernel.
The {\it Poisson envelope} of $A$ is 
$$
PA:=\frac{SLA}{<<IA>>}
$$
The inclusion $A=SA/IA\subset PA$ has the following universal property. If
$P$ is a Poisson algebra and $f:A@>>>P$ is a homomorphism of commutative 
algebras, then there is a unique Poisson homomorphism $PA@>>>P$ which extends
$f$. In other words $A\mapsto PA$ is left adjoint to the forgetful functor
$\poiss@>>>\comm$ from Poisson to commutative algebras. 
One checks that if $A=SV/I$ is any presentation of $A$ as a quotient of a 
symmetric algebra, then $SLV/<<I>>$ has the same universal property as and 
is therefore isomorphic to $PA$. In particular 
$$
PSV=\fpoiss V
$$
It follows from \thetag{11} that if $I\subset SV$ is as above then
$<<I>>\subset SLV$ is homogeneous for the grading \thetag{7}, whence 
$PA$ inherits a grading:
$$
PA=\bigoplus_{n\ge 0}P_nA
$$
For example 
$$
P_nSV=S_nLV\qquad (n\ge 0)\tag{12}
$$
In particular
$$
P_1SV=SV\otimes L_1V=SV\otimes\Lambda^2V=\Omega^2_{SV}\tag{13}
$$
is the module of $2$-differential forms. If follows from \thetag{13} and 
\thetag{11} that for every commutative algebra $A$,
$$
\split
P_1A=\frac{S_1(LA)}{IAS_1(LA)+<\{A,IA\}>}\\
    =\frac{\Omega^2_{SA}}{IA\Omega^2_{SA}+\Omega^1_A\wedge dIA}
    =\Omega^2_A
\endsplit
$$
Under this isomorphism,
$$
\{a,b\}=da\wedge db\in P_1A
$$
For another interpretation of $P_1A$ consider the analogy $L^A$ of the functor
$L$ for $A$-modules and $A$-Lie algebras. If $M$ is an $A$-module, then
$L^AM$ carries a grading defined exactly as in \thetag{6}. We have $L^A_0M=M$,
$L^A_1M=\Lambda^2M$, and in particular
$$
L_1^A\Omega^1_A=\Omega^2_A=P_1A\tag{14}
$$
Theorem 4.1.3 of [2] says that a generalization of \thetag{14} holds
for smooth algebras. Namely it is asserted that for every $n\ge 0$,
$P_nA$ is isomorphic as an $A$-module
to the homogeneous part of degree $n$ of the symmetric $A$-algebra on 
$L^A\Omega^1_A$:
$$
P_nA\cong S_n^AL_+^A\Omega^1_A\qquad(n\ge 0)\tag{15}
$$
However the proof of this assertion in [2] has a gap, as the isomorphism
given there is not well-defined. Indeed the map in question sends the
element 
$$
P_nA\owns b\{a_0,\{a_1,\dots,\{a_{n-1},a_n\}\dots\}\}\quad(a_i\in A)
$$
to the element
$$
b[da_0,[da_1,\dots[da_{n-1},da_n]\dots]]\in S_n^AL_+^A\Omega^1_A
$$
However a calculation shows that this rule maps 
$$
0=\{a_1,a_3\{a_2,a_4\}\}+\{a_1,a_2\{a_3,a_4\}\}-\{a_1,\{a_2a_3,a_4\}\}
$$
to the element
$$
[da_1,da_3][da_2,da_4]+[da_1,da_2][da_3,da_4]
$$
which is nonzero in general. I do not know whether the isomorphism 
\thetag{15} 
still holds for every smooth
algebra $A$. It certainly holds for symmetric algebras, as is immediate from
\thetag{12}. We show in Theorem 1.4 below that it also holds for local 
smooth algebras. For a weaker version of \thetag{15} which holds globally, 
see Theorem 1.6. The following lemma is well-known.
\bigskip
\proclaim{Lemma 1.3} Let $X$ be a set, $V=\oplus_{x\in X}kx$ the free 
vectorspace on $X$. Then the set 
$$
Y:=X\cup\bigcup_{n=1}^\infty\{[x_1,[x_2,[\dots ,[x_n,x_{n+1}]\dots ]]]: x_i\in 
X\}
$$
generates $LV$ as a vectorspace. In particular there is a basis $Z$ of $LV$
such that $X\subset Z\subset Y$.
\endproclaim
\demo{Proof} Straightforward induction.\qed\enddemo
\bigskip
\proclaim{Theorem 1.4} Let $A$ be a local smooth algebra. Then $A$ satisfies
\thetag{15}. 
\endproclaim
\demo{Proof}Let $x_1,\dots,x_n\in A$  be a regular system of parameters
and $V=\oplus_{i=1}^nkdx_i\subset\Omega^1_A$. We have
$$
S^AL^A_+\Omega^1_A=A\otimes SL_+V
$$
Put $L=LV$; then
$$
\Omega^1_A=A\otimes V,\ \  \Omega^1_{S{L}_+}=SL_+\otimes L_+, \ \
\Omega^1_{A\otimes S{L}_+}=A\otimes SL_+\otimes L\tag{16}
$$
Consider the permutation isomorphism
$$
\tau:\Omega^1_A\otimes SL_+=A\otimes V\otimes SL_+\cong A\otimes 
SL_+\otimes V
$$
Under the identifications \thetag{16} the de Rham derivation 
$d_{A\otimes S{L}_+}$ is identified with 
$$
D:=\tau\circ(d_A\otimes id_{S{L}_+})+id_A\otimes d_{S{L}_+}
$$
Put 
$$
\phi=id_{A\otimes S{L}_+}\otimes [,]:A\otimes SL_+\otimes \Lambda^2L
@>>>A\otimes SL_+
$$
One checks that 
$$
\{p,q\}:=\phi(Dp\wedge Dq)
$$
is a Poisson bracket. Note that for 
$1\le i_1,\dots,i_{r+1}\le n$, we have
$$
\{x_{{i}_1},\{x_{{i}_2},\dots,\{x_{{i}_r},
x_{{i}_{r+1}}\}\dots\}\}=
[dx_{{i}_1},[dx_{{i}_2},\dots,
[dx_{{i}_r},dx_{{i}_{r+1}}]\dots]]
$$
It suffices to show that the Poisson algebra $(A\otimes SL_+,\{,\})$
together with the inclusion 
$A=A\otimes k=A\otimes S_0L_+\subset A\otimes SL_+$ has the universal
property of $PA$. Let $P$ be a Poisson algebra and $f:A@>>>P$ a 
homomorphism of commutative algebras. Write $p_i=f(x_i)$. By lemma 1.3, we
may extend $B_0=\{dx_1,\dots,dx_n\}$ to a homogeneous basis $B$ of $L$ such
that every element of $B':=B\backslash B_0$ be of the form 
$[dx_{{i}_1},[dx_{{i}_2},\dots, [dx_{{i}_r},dx_{{i}_{r+1}}]]\dots]$ ($r\ge 1$).
View $P$ as an $A$-module via $f$ and consider the $A$-module homomorphism
$\theta:L^A_+\Omega^1_A=A\otimes L_+@>>>P$ defined on elements of $B'$ by 
$$
\theta[dx_{{i}_1},[dx_{{i}_2},\dots,
[dx_{{i}_r},dx_{{i}_{r+1}}]\dots]]=\{p_{{i}_1},\{p_{{i}_2},\dots,
\{p_{{i}_r},p_{{i}_{r+1}}\}\dots\}\}\tag{17}
$$
Note that, as defined,
$$
\theta[l_1,l_2]=\{\theta l_1,\theta l_2\}\qquad (l_1,l_2\in B')\tag{18}
$$
Indeed by lemma 1.3, the two sides of this identity are defined by the same
linear combination of the elements of $B'$ and of their images. The 
prescription \thetag{17} together with the prescription that $\theta$ 
extend $f$,
determine a unique map $\theta:A\otimes SL_+@>>>P$ which satisfies 
\thetag{18} for $l_1,l_2\in B'':=\{x_1,\dots,x_n\}\cup B'$. It follows
that $\theta$ is a Poisson homomorphism; uniqueness is clear. \qed\enddemo 
\bigskip
\subhead{1.5. A weaker version of property \thetag{15}}\endsubhead
Let $V$ be a vectorspace. Combining the two gradings \thetag{7}, 
\thetag{8} we obtain a bigrading
$$
SLV=\bigoplus_{p,q\ge 0}S_q^pLV\tag{19}
$$
where 
$$
S_q^pLV:=S^pLV\cap S_qLV
$$
Now let $A$ be a commutative algebra, and consider the projection
$$
\pi:SLA\fib PA\tag{20}
$$
The map $\pi$ is homogeneous for the $||_*$-degree but not for the 
$||^*$-degree. However the ideal
$$
\Cal H^n:=\pi(\bigoplus_{p\ge n}S^pLA)\qquad (n\ge 0)
$$
is homogeneous with respect to $||_*$ and therefore the graded
ring $GPA=G_{\Cal H}PA$ is actually bigraded
$$
GPA=\bigoplus_{p,q\ge 0}G^p_qPA
$$
On the other hand \thetag{19} carries over to the free Lie algebra of any
$A$-bimodule. In particular $S^AL^A\Omega^1_A$ is a bigraded algebra.
\bigskip
\proclaim{Theorem 1.6} Let $A$ be a commutative algebra. Then for the bigraded
structures defined in 1.5, there is a natural
surjective homomorphism of bigraded algebras
$$
\phi:S^AL^A\Omega^1_A\fib GPA
$$
If furthermore $A$ is smooth, then $\phi$ is an isomorphism.
\endproclaim
\demo{Proof} The map \thetag{20} is the quotient by $<<IA>>$, whence it
factors through a map
$$
\overline{\pi}:A\otimes SL_+A=\frac{SLA}{<IA>}\fib PA
$$
In particular $\overline{\pi}$ induces a homogeneous, surjective homomorphism
of graded $A$-modules
$$
p:A\otimes L_+A@>\overline{\pi}>>\Cal H^1\fib G_*^1PA\
$$
Write $\rho:A@>>>SA$ for the canonical inclusion. The ideal $IA$ is generated
by the elements
$$
u(a,b):=\rho(ab)-\rho a\rho b\qquad (a,b\in A)
$$
Let $h_i(a_1,\dots,a_n;b,c)$ be the homogeneous part of 
$||^*$-degree one of the element $g_i(\rho a_1,\dots,\rho a_n; w(b,c))$ of
\thetag{10}. By
\thetag{11}, the elements $h_i(a_1,\dots,a_n;b,c)$, 
$1\le i\le n$, $a_i,b,c\in A$ generate $\ker p$ as an $A$-module.
A calculation shows that 
$$
\gather
h_i(a_1,\dots,a_n;b,c)=\\
1\otimes\{\rho a_1,\dots,\{\rho a_i,\{\rho(bc),
\{\rho a_{i+1},\dots,\{\rho a_{n-1},\rho a_n\}\dots\}\}\}\dots\}-\\
\rho b\otimes \{\rho a_1,\dots,\{\rho a_i,\{\rho c,
\{\rho a_{i+1},\dots,\{\rho a_{n-1},\rho a_n\}\dots\}\}\}\dots\}-\\
\rho c\otimes \{\rho a_1,\dots,\{\rho a_i,\{\rho b,
\{\rho a_{i+1},\dots,\{\rho a_{n-1},\rho a_n\}\dots\}\}\}\dots\}\quad
\endgather
$$
It follows from this that 
$$
M:=\ker (A\otimes A\fib \Omega^1_A)\oplus\ker p
$$
is the Lie ideal generated by $\ker(A\otimes A\fib\Omega^1_A)$ in the $A$-Lie 
algebra $A\otimes LA$ . Hence $(A\otimes LA)/M=L^A\Omega^1_A$, and
$$
 L^A_+\Omega^1_A\cong G^1PA\tag{21}
$$
The map $\phi$ of the theorem is that induced by \thetag{21}; it is
surjective because $G^1PA$ generates $GPA$ as an $A$-algebra. Assume now
that $A$ is smooth; we must prove that $\phi$ is injective. This is a local
question, so we may further assume that $A$ is local. Let $x_1,\dots,x_n$
be a regular system of parameters. By the proof of 1.4 there is an isomorphism
$\psi:PA\iso S^AL^A_+\Omega^1_A$ such that
$$
\psi\{x_{{i}_1},\{x_{{i}_2},\dots,
\{x_{{i}_r},x_{{i}_{r+1}}\}\dots\}\}=[dx_{{i}_1},[dx_{{i}_2},\dots,
[dx_{{i}_r},dx_{{i}_{r+1}}]\dots]]
\tag{22}
$$
Furthermore the induced map $G\psi:GPA@>>>S^AL^A_+\Omega^1_A$ still verifies 
\thetag{22}. Thus $G\psi\circ\phi$ is the identity map, because it is
so on generators. In particular, $\phi$ is injective.\qed
\enddemo
\bigskip
\head{2. Local quantization of the Poisson envelope}\endhead

\subhead{2.0. PBW quantization}\endsubhead
Let $\g$ be a Lie algebra, $S\g$ and $U\g$ the symmetric and universal 
enveloping algebras, and consider the symmetrization map
$$
e:S\g@>>>U\g\quad e(g_1\dots g_n)=\frac{1}{n!}\sum_{\sigma\in\frak S_n}
g_{\sigma(1)}\dots g_{\sigma(n)}\tag{23}
$$
By the Poincar\'e-Birkhoff-Witt theorem, the associative product
$$
B:S\g\otimes S\g@>>>S\g,\quad B(x\otimes y)=e^{-1}(exey)
$$
decomposes as a sum
$$
B=\sum_{p=0}^\infty B_p\ \ \text{where\ \ } B_p(S^n\g)\subset S^{n-p}\g
\ \ (n,p\ge 0)\tag{24}
$$
We have $B_0(x\otimes y)=xy$, $B_1(x\otimes y)=\frac{1}{2}\{x,y\}$. Explicit
formulas for all the $B_p$ are given in [1]. It also proved in {\it loc. cit.}
that for each $p\ge 0$, $B_p$ is a differential operator of order $\le p$. We
call the map
$$
B(\hbar):=\sum_{n\ge 0}B_n{\hbar}^n:S\g\otimes S\g [[\hbar]]@>>>
S\g[[\hbar]]
$$
the PBW {\it quantization}. The next lemma establishes the properties of the
product $B$ with respect to the commutator filtration \thetag{1} in the
case when $\g$ is free.
\bigskip
\proclaim{Lemma 2.1} Let $V$ be a vectorspace, $LV$ the free Lie algebra,
$TV$ the tensor algebra, and $e$ as in \thetag{23}. Then
\roster
\item $e(\fpoiss_{\ge n}V)=F_nTV$
\item The operator $B_p:\fpoiss V\otimes \fpoiss V@>>>\fpoiss V$ of 
\thetag{24} is 
homogeneous of degree $+p$ for the $||_*$-degree \thetag{7}.
\endroster
\endproclaim
\demo{Proof} One checks that if $\g$ is any Lie algebra and $U\g$ its 
enveloping
algebra, then for $F_0\g=\g$, $F_d\g=[\g,F_{d-1}\g]$ ($d\ge 1$)
we have
$$
F_nU\g=\sum_r\sum_{{d}_1+\dots +{d}_r\ge n}F_{{d}_1}\g\cdot\dots\cdot F_{{d}_r}\g
$$
For $\g=LV$, we obtain
$$
F_nTV=\sum_r\sum_{{d}_1+\dots +{d}_r\ge n}L_{{d}_1}V\cdot\dots\cdot L_{{d}_r}V
\tag{25}
$$
From \thetag{25} and \thetag{9} it is clear that 
$e(\fpoiss_{\ge n}V)\subset F_nTV$.
We must show that $F_nTV\subset e(\fpoiss_{\ge n}V)$. Consider the following 
subspace of $F^nTV$
$$
\Cal A_{p,n}:=\sum_{r\le p}\sum_{{d}_1+\dots +{d}_r\ge n}L_{{d}_1}V\cdot\dots
\cdot L_{{d}_r}V
$$
Clearly $F_nTV=\cup_{p\ge 1}\Cal A_{p,n}$ An inductive argument similar to
that of the usual proof of the surjectivity of $e$ shows that for $p\ge 1$,
$\Cal A_{p,n}\subset e(\fpoiss_{\ge n} V)$. This proves assertion (1). 
Assertion (2) follows by counting degrees in the formula for $B_p$ given 
in [1] 1.1.\qed
\enddemo
\proclaim{Corollary 2.2}(Compare [2], 3.4.7) The natural map 
$$\fpoiss V=PSV\iso\bigoplus_{n=0}^\infty F_nTV/F_{n+1}TV$$
is an isomorphism.\qed
\endproclaim
\bigskip
\subhead{2.3. PBW quantization of $PA$ for $A$ local and smooth}\endsubhead
Let $(A,\Cal M)$ be a smooth local commutative algebra and
$X=\{x_1,\dots,x_n\}\subset\Cal M$ a regular system of parameters. Set 
$V=k^n$. We are going to combine Theorem 1.4 and the PBW 
quantization of $\fpoiss V$ to obtain an associative product
$$
B^X(\hbar)=\sum_{p=0}^\infty B^X_p\hbar^p:PA\otimes PA[[\hbar]]@>>>PA[[\hbar]]
$$
Because the map $B_p:\fpoiss V\otimes \fpoiss V@>>>\fpoiss V$ is a 
differential operator, it is continuous with respect to the topology of any
ideal $I\subset \fpoiss V$. Applying this for $I=V\cdot\fpoiss V$ 
and completing we obtain
the horizontal solid arrow in the following commutative diagram
$$
\xymatrix{
PA\otimes PA\ar[d]_{{\iota}_X\hat{\otimes}{\iota}_X}\ar@{..>}[r]^{{B}^{X}_{p}}
 &PA\ar[d]^{{\iota}_X}\\
k[[t]]\hat{\otimes}SL_+V\hat{\otimes} k[[t]]\hat{\otimes}SL_+V
\ar[r]_(.56){\quad \hat{B_p}} & k[[t]]\hat{\otimes} SL_+V
}
\tag{26}
$$    
Here $\hat{\otimes}$ is the completed tensor product and $k[[t]]$ is shorthand 
for $k[[t_1,\dots,t_n]]$. The map $\iota_X$ is 
the composite
$$
\iota_X:PA\overset{\alpha_X}\to\iso A\otimes SL_+V\hookrightarrow \hat{A}
\hat{\otimes} SL_+V\overset{j_X\otimes 1}\to\iso k[[t]]\hat{\otimes}SL_+V
$$
where $\alpha_X$ is the isomorphism of 1.4, $\hookrightarrow$ is the passage 
to completion and $j_X:\hat{A}\cong k[[t]]$ is the isomorphism determined by $x_i\mapsto t_i$ ($i=1,\dots,n$). Because $A\hookrightarrow \hat{A}$ is 
injective, so are both
vertical maps in \thetag{26}. The map $B^X_p$ is defined by the following
lemma.
\bigskip
\proclaim{Lemma 2.4} The map $\hat{B}_p$ of \thetag{26} sends the image
of $\iota_X\hat{\otimes}\iota_X$ to the image of $\iota_X$.
\endproclaim
\demo{Proof} Let $Z=X\cup\{y_1,y_2,\dots\}$ be a basis as that of Lemma 1.3. 
Every monomial
on the elements of $Z$ can be written as $x^{\alpha'}y^{\alpha''}=
\prod_{i=1}^nx_i^{\alpha'(i)}\prod_{j=1}^\infty y_j^{\alpha''(j)}$ for
some multi-indices $\alpha':\{1,\dots,n\}@>>>\zahl_{\ge 0}$ and 
$\alpha'':\zahl_{\ge 1}@>>>\zahl_{\ge 0}$ with $\alpha''(n)=0$ for $n>>0$.
Let $\frac{\partial^{|\alpha|}}
{\partial x^{\alpha'}\partial y^{\alpha''}}$ be the higher derivation with
symbol $x^{\alpha'}y^{\alpha''}$.
Because $B_p:SLV\otimes SLV@>>>SLV$
is a bidifferential operator of order $\le p$, it can be written as an 
$SLV$-linear combination of cup products of higher derivations with respect
to the basis $Z$
$$
B_p=\sum_{|\alpha|,|\beta|\le p}c_{\alpha,\beta}
\frac{\partial^{|\alpha|}}{\partial x^{\alpha'}\partial 
y^{\alpha''}}\cup\frac{\partial^{|\beta|}}
{\partial x^{\beta'}\partial y^{\beta''}}
$$
Since each of the higher derivations above maps $A\subset k[[t]]$ to itself,
so does $\hat{B}_p$.\qed
\enddemo
\bigskip
\proclaim{Theorem 2.5} Let $(A,\Cal M)$ be a smooth local algebra,
$X\subset \Cal M\cdot A_{\Cal M}$ a regular system of parameters, $p\ge 0$,
$B^X_p$ as in \thetag{26} and $\widehat{PA}=\varprojlim_dPA/P_{>d}A=
\prod_{d=0}^\infty P_dA$. Then 
\roster
\item 
$B^X(\hbar)=\sum_{p=0}^\infty B^X_p\hbar^p:PA\otimes PA[[\hbar]]@>>>
PA[[\hbar]]$ 
is associative.
\item $B^X_p$ is a differential operator of order $\le p$.
\item $B^X_p(P_nA\otimes P_mA)\subset P_{n+m+p}A$.
\item The map $B^X(\hbar)$ induces a continuous associative product
$$B^X(1):=\sum_{p=0}^\infty B^X_p:\widehat{PA}\hat\otimes\widehat{PA}
@>>>\widehat{PA}$$
\item For the associative algebra $Q_X=(\widehat{PA},B^X(1))$, we have 
$$F_nQ_X=\prod_{d\ge n}P_dA$$
\item There is a commutative diagram of monomorphisms

$$
\xymatrix{
k\{t_1,\dots,t_n\}\ar[r]\ar[dr]_{inc} & Q_X\ar[d]\\
&k\{\{t_1,\dots,t_n\}\}
}
$$ 
where $inc$ is the canonical inclusion of the noncommutative polynomials
into the noncommutative power series.
\endroster
\endproclaim
\demo{Proof} Assertions (1), (2) and (3) are immediate from the analogous
properties of the PBW quantization; (4) follows from (3), and (5) from Lemma 
2.1. It is clear 
from the definition of $Q_X$ that there is a diagram as that in (5) but with 
$\widehat{SLV}:=(\prod_{n\ge 0}k[[t]]\hat{\otimes}SL_+V,\hat{B})$ substituted 
for 
$k\{\{t\}\}$. Note that, for $I=\sum t_i\cdot k[t_1,\dots, t_n]$,
$\widehat{SLV}$ is the completion of $SLV$ with 
respect to the filtration $\{ I\cdot SLV+S_{\ge n}LV:n\ge 0\}$, and that
$J:=e^{-1}(<t_1,\dots,t_n>)=I\oplus S_+LV$. By lemmas 2.1 and 2.6
$\widehat{SLV}\cong\varprojlim_{n}TV/e(J)^n=k\{\{t_1,\dots,t_n\}\}$.\qed
\enddemo

\proclaim{Lemma 2.6} Let $G =\oplus_{n=0}^\infty G_n$ be a graded
commutative algebra. Assume $G$ is additionally equipped with an
associative --but not necesarily commutative- product 
$$
\Phi=\sum_{p=0}^\infty \Phi_p:G\otimes G@>>>G
$$
such that $\Phi_0$ is the original commutative product, and that for each 
$p\ge 1$, $\Phi_p$ is a bidifferential operator. Let $I\subset G_0$ be an 
ideal for $\Phi_0$, and $J\subset G$ the $\Phi$-ideal it 
generates. Then the linear topologies induced on 
the underlying vectorspace of $G$ by the filtrations 
$\{J^n+G_{\ge n}:n\ge 0\}$ and 
$\{I^nG+G_{\ge n}:n\ge 0\}$ coincide. 
\endproclaim
\demo{Proof} 
It suffices to prove that for each $d\ge 0$ the filtrations 
$\{J^n+G_{\ge d+1}/G_{\ge d+1}: n\ge 0\}$ and 
$\{I^nG+G_{\ge d+1}/G_{\ge d+1}:n \ge 0\}$ of $G/G_{\ge d+1}$ are equivalent. 
Thus we may assume that $G_m=0$ for $m\ge d+1$. Hence $\Phi$ is a 
bidifferential operator. Let $\alpha$ be the order of $\Phi$. We write
$x\star y:=\Phi(x\otimes y)$, and if $X\subset G$ is any subspace, we put
$X^{\star n}$ for the subspace generated by all products $x_1\star\dots\star x_n$
with $x_i\in X$. Let $i\in I$, $p\ge 1$, $n,r\ge 0$, 
$F_i(x):=i\star x$. Because $F_i$ is a differential operator of order $\le\alpha$,
$$
F_i(I^{p\alpha+r}G_n)\subset I^{p\alpha+r+1}G_n+I^pG_{\ge n+1}\tag{27}
$$
Using \thetag{27} one checks by induction that for 
$(c_r\dots c_0):=\sum_{i=0}^rc_i\alpha^i$,
$$
I^{\star (c_r\dots c_0)}\subset M_{({c}_r\dots {c}_0)}:=\sum_{j=0}^rI^{({c}_r
\dots {c}_j)}G_{\ge j}
+G_{\ge r+1}\tag{28}
$$
We remark that $M_{({c}_r\dots {c}_0)}$ is an ideal for both $\star$ and the 
original product. Hence for $N\ge ({c}_r\dots {c}_0)+d$ and $r\ge d$,
$$
J^{\star N}\subset (I\oplus G_{\ge 1})^{\star N}\subset 
<I^{\star ({c}_r\dots {c}_0)}>_\star
\subset M_{({c}_r\dots {c}_0)}\subset I^{(c_r\dots c_d)}G
$$
where the subindex $_\star$ denotes two sided ideal generated by the product 
$\star$. Now using \thetag{28}, and noting that 
$I^{\star n}\star G_d=I^nG_d$
and that in general the projection $G\fib G_j$ maps $I^{\star n}\star G_j$ 
surjectively onto $I^nG_j$, one checks that, for $r\ge d$ 
$$
M_{(c_r\dots c_0)}=\sum_{j=0}^d I^{\star (c_r\dots c_j)}\star G_{d-j}
$$
It follows that
$$
J^{\star n}\supset M_{n{\alpha}^d}\supset I^{n{\alpha}^d}\qed
$$ 
\enddemo

\bigskip
\head{3. Smooth nilcommutative and nil-Poisson algebras}\endhead

\subhead{3.0 Nil-Poisson algebras}\endsubhead
Let $P$ be a Poisson algebra. Put $F_0P=P$ and inductively
$$
F_{n+1}P:=\sum_{i=1}^nF_iPF_{n+1-i}P+\sum_{i=0}^n<\{F_iP,F_{n-i}P\}>\tag{29}
$$
for $n\ge 0$. This is the Poisson analogue of the commutator filtration 
\thetag{1}. For example if $A$ is any commutative algebra then
$$
F_rPA=P_{\ge r}A=\bigoplus_{n=r}^\infty P_nA\qquad (r\ge 0)
$$
The analogue of a nilcommutative algebra of order $\le d$
is called a {\it nil-Poisson} algebra of order $\le d$. The category of
ni-Poisson algebras of order $\le d$ is $NP_d$. We put 
$NP=\cup_{d\ge 0}NP_d$. Formal $d$-smoothness and
$d$-smoothness for objects of $NP_d$ are the obvious analogues of the same
properties for objects of $NC_d$ as defined in \S1.
If $A$ is (formally) smooth in the commutative sense, then $PA/P_{>d}A$ is
(formally) $d$-smooth. It turns out that every (formally) $d$-smooth Poisson 
algebra is of this form;
see proposition 3.3. The following Lemma is the analogue of [2], 1.2.7 for
Poisson algebras.
\bigskip
\proclaim{Lemma 3.1} Let $P\in NP_d$ and $f:P@>>>P$ a Poisson endomorphism.
Assume the induced map $P/F_1P@>>>P/F_1P$ is the identity. Then the restriction
of $f$ to $F_{d+1}P$ is the identity also.
\endproclaim
\demo{Proof} Consider the map $D:P@>>>P$, $Dp:=fp-p$. We have
$$
\gather
D\{p,q\}=\{Dp,q\}+\{p,Dq\}-\{Dp,Dq\}\tag{30}\\
D(pq)=pDq+qDp-DpDq
\endgather
$$
By hypothesis, $DP\subset F_1P$; it follows from this, using \thetag{30},
\thetag{29} and induction, that for $n\ge 0$, $D(F_nP)\subset F_{n+1}P$.
In particular $D(F_dP)=0$.\qed
\enddemo
\bigskip
\remark{Remark 3.2} The proof of the lemma above still applies if 
one substitutes $NC_d$ for $NP_d$ and ``algebra endomorphism'' for ``Poisson 
endomorphism''. This gives an alternate proof of [2], 1.2.7.
\endremark
\bigskip
\proclaim{Proposition 3.3} Let $P\in NP_d$, $A=P/F_1P$. Then the following
conditions are equivalent
\item{i)} $P$ is $d$-formally smooth.
\smallskip
\item{ii)} $A$ is $0$-formally smooth and $PA/P_{>d}A\cong P$.
\smallskip
The same holds if we replace ``formally smooth'' by ``smooth'' in both i) and 
ii).
\endproclaim
\demo{Proof} That ii)$\Rightarrow$i) is clear, as is that i) implies $A$
is formally smooth. Use the formal smoothness of $A$ to obtain a section
$s:A@>>>P\in\comm$ of the projection $P\fib A$, and then the universal property
of $PA$ to lift $s$ to a map of extensions $\alpha:PA/P_{>d}A@>>>P$. To prove 
that if i) holds then $\alpha:PA/P_{>d}A@>>>P$ 
is an isomorphism, note that, by the hypothesis on $P$, there is
a map $\beta:P@>>>PA/P_{>d}A$ which descends to the identity of $A$. One is 
thus reduced to showing that $\alpha\beta$ and $\beta\alpha$ are isomorphisms.
This follows from lemma 3.1.
\qed
\enddemo
\bigskip
Part (1)$\Rightarrow$(2) of the following theorem is due to Kapranov
([2] 4.2.1); we give a new proof.

\bigskip
\proclaim{Theorem 3.4} Let $R\in NC_d$, $G=\oplus_{n=0}^dG_n$ the associated
graded Poisson algebra, $A=G_0$, $\pi:R@>>>A$ the projection. 
The following conditions are equivalent
\roster
\item $R$ is $d$-smooth.
\item $A$ is $0$-smooth and the canonical map $PA/P_{>d}A@>>>G$ is an 
isomorphism.
\item For every maximal ideal $\Cal M\subset A$, there is a regular system
of parameters $X\subset\Cal M\cdot A_{\Cal M}$ such that for $Q_X$ as
in theorem 2.5, the \O re localization of
$R$ at $\pi^{-1}(\Cal M)$ is isomorphic to $Q_X/F_{>d}Q_X$.
\endroster
\endproclaim
\demo{Proof} Assume (1) holds. Then clearly $A$ is $0$-smooth. Let 
$\Cal M\subset A$ be a maximal ideal, $\bold{M}=\pi^{-1}(\Cal M)$,
$X\subset\Cal M\cdot A_{\Cal M}$
a regular system of parameters and $Q=Q_X/F_{>d}Q_X$. Because the \O re 
localization
$R_{\bold{M}}$ is $d$-smooth, the identity of $A_{\Cal M}$ can be lifted
to a map $f:R_{\bold{M}}@>>>Q$. By lemma 3.1, the map induced by 
$f$ at the graded level is an isomorphism, whence $f$ is an isomorphism.
We have just proved that (1)$\Rightarrow$(3). 
It is clear that (3)$\Rightarrow$(2). We prove next that (2)$\Rightarrow$(1),
by induction on $d\ge 0$. The case $d=0$ is tautological.
Assume $d\ge 1$ and that the theorem is true for $d-1$. Let $R\in NC_d$ satisfy
the hypothesis of the theorem. Then by inductive assumption 
$R_{d-1}=R/F_{d}R$ is $d-1$-smooth. Let $\pi:U\fib R_{d-1}$ be the universal central
extension as defined in [2] 1.3.6. By [2] 1.6.2, $U$ is $d$-smooth. One
checks that $\ker\pi=F_dU$. By [2] 1.3.8, there is a map 
$\alpha:U@>>>R$ which induces the identity of $R_{d-1}$. Consider the 
Poisson homomorphism $\beta$ induced by $\alpha$ at the associated graded
level. Because (1)$\Rightarrow$(2), $\beta$ is an endomorphism of 
$PA/P_{>d}A$. By virtue of Lemma 3.1, because
$\beta$ induces the identity modulo $P_dA$ and is homogeneous,  
it has to be the identity. Thus $\alpha$ is an isomorphism.\qed
\enddemo


\bigskip
\Refs

\ref\no{1}\by Corti\~nas, G.\paper An explicit formula for PBW quantization
\paperinfo math.QA/0001127
\endref

\ref\no{2}\by Kapranov, M.\paper Noncommutative geometry based on 
commutator expansions\jour J. reine angew. Math.\vol 505\yr 1998
\pages 73-118
\endref

\ref\no{3}\by Serre, J.P.\book Lie algebras and Lie groups\publ W. A. 
Benjamin, Inc.\yr 1965
\endref


\endRefs
\enddocument